\documentclass[11pt,a4paper]{article}
\usepackage{amsfonts, amssymb, amsmath, fancyhdr, url, amsthm}
\usepackage{enumerate}
\usepackage{xcolor}
\usepackage{geometry}
\usepackage{tocloft}
\usepackage[backref=page]{hyperref}
\renewcommand*{\backref}[1]{}
\renewcommand*{\backrefalt}[4]{%
	\ifcase #1 (Not cited.)%
	\or        (Cited on page~#2.)%
	\else      (Cited on pages~#2.)%
	\fi}

\hypersetup{
	colorlinks=true,
	linkcolor=blue,
	filecolor=blue,
	citecolor = red,
	urlcolor=cyan,
}

\newcommand{\al}{\alpha}

\newcommand{\f}{\varphi}

\newcommand{\Lie}{\operatorname{Lie}}

\newcommand{\K}{K\"ahler}

\newcommand{\eg}{{\bf e.g. }}

\newcommand{\CC}{\mathbb{C}}

\newcommand{\RR}{\mathbb{R}}

\def\eqref#1{(\ref{#1})}

\newcommand{\C}{{\mathbb C}}
\newcommand{\R}{{\mathbb R}}

\newcommand{\6}{\partial}
\def\1{\sqrt{-1}\:}

\newcommand{\cntrct}                
{\hspace{2pt}\raisebox{1pt}{\text{$\lrcorner$}}\hspace{2pt}}

\newcommand{\ii}{\mathrm{i}}

\renewcommand{\bar}{\overline}
\renewcommand{\phi}{\varphi}
\renewcommand{\epsilon}{\varepsilon}

\newcommand{\Id}{\operatorname{Id}}

\numberwithin{equation}{section}

\makeatletter

\@addtoreset{equation}{section}
\@addtoreset{footnote}{section}
\makeatother

\def\blacksquare{\hbox{\vrule width 5pt height 5pt depth 0pt}}
\def\endproof{\blacksquare}

\newcommand{\samethanks}[1][\value{footnote}]{\footnotemark[#1]}

\newcounter{Mycounter}[section]
\newcounter{lemma}[section]
\newcounter{claim}[section]

\newcounter{sublemma}[section]

\newcounter{corollary}[section]

\newcounter{theorem}[section]

\newcounter{conjecture}[section]

\newcounter{proposition}[section]

\newcounter{definition}[section]

\renewcommand{\thedefinition} {{Definition~\thesection.\arabic{definition}}}

\newcommand{\definition}{%

     \setcounter{definition}{\value{Mycounter}}

     \refstepcounter{definition}

     \stepcounter{Mycounter}

     {\vspace{.1in} \noindent \bf \thedefinition:\ }}

\newcounter{example}[section]

\renewcommand{\theexample}{{Example \thesection.\arabic{example}}}

\newcommand{\example}{%

     \setcounter{example}{\value{Mycounter}}

     \refstepcounter{example}

     \stepcounter{Mycounter}

     {\vspace{.1in}\noindent \bf \theexample:\ }}

\newcounter{remark}[section]

\renewcommand{\theremark}{{Remark \thesection.\arabic{remark}}}

\newcommand{\remark}{%

     \setcounter{remark}{\value{Mycounter}}

     \refstepcounter{remark}

     \stepcounter{Mycounter}

     {\vspace{.1in}\noindent \bf \theremark:\ }}

\newcounter{problem}[section]

\newcounter{question}[section]

\begin{document}

\begin{center}
{\LARGE\bf
Deformations of Vaisman manifolds 
}\\[3mm]
{\large\bf
L. Ornea\footnote{Partially supported by by Romanian Ministry of Education and Research, Program PN-III, Project number PN-III-P4-ID-PCE-2020-0025, Contract  30/04.02.2021} and V. Slesar\samethanks{}

}
\end{center}

\hfill

{\small
\hspace{0.15\linewidth}
\begin{minipage}[t]{0.7\linewidth}
{\bf Abstract}
We construct a type of transverse deformations of a Vaisman manifold, which preserves the canonical foliation. For this construction we only need a basic 1-form with certain properties. We show that such basic 1-forms exist in abundance.  \\

\end{minipage}
}

\hspace{5mm}



\tableofcontents

\date{\today }

\section{Introduction}

A { Vaisman manifold} is a particular Hermitian manifold $(M, J, g)$ with its fundamental form $\omega(\cdot, \cdot):=g(\cdot, J\cdot)$ satisfying the relation $d\omega=\theta\wedge \omega$, for a nonzero one-form $\theta$, which is parallel with respect to the Levi-Civita connection of the metric $g$.  The one-form $\theta$ is called the {\bf Lee form}. We assume throughout the paper that the manifold $M$ is closed, connected, with $\dim_\C M\ge 2$.

Since a parallel form is closed,  Vaisman manifolds are locally conformally {\K} (LCK). Note that many of the known (LCK) manifolds are in fact Vaisman. For the main properties and  examples of LCK manifolds,  see e.g. \cite{Dr-Or}, \cite{ov_lee}.

Vaisman manifolds bear a holomorphic foliation of complex dimension 1, generated by the Lee and anti-Lee fields $\theta^\sharp$ and $J\theta^\sharp$, usually called { the canonical foliation}. It is locally Euclidean and transversally \K.

Unlike \K\ structures, LCK structures are not stable to small deformations (\cite{Be}). However, Vaisman structures are preserved by some particular type of deformations obtained by perturbing the Lee vector field in the Lie algebra of the closure of the group generated by its flow, see \cite{Or-Ve1}.

Other types of deformations of Vaisman structures are obtained by fixing the complex structure and the Lee vector field and changing the Lee form (\cite{Be}) (these are related to deformations of Sasaki structures, see also \cite{Bo-Ga}).

In this paper, we adopt another strategy. We construct a deformation $(M,J_{t},g_{t})$  of the Vaisman structure, with $J_{0}=J$, $g_{0}=g$ and $t$ sufficiently small, in such a way that the canonical foliation is not affected and the deformation concerns only the transverse orthogonal complement and the transverse K\"ahler geometry. This type of deformation is related  to the deformations of second type on Sasakian manifolds, as defined in \cite{Be}. For our construvtion, we need a basic 1-form with certain natural properties and we indicate a way of producing such 1-forms (Subsection \ref{Construction}). This procedure does not cover the whole space of deformations, still it gives a way of obtaining new Vaisman structures out of a given one (see \ref{g_t_vai}).  We end with a completely worked-out application of the method to a classical Hopf surface $\frac{\C^2\backslash 0}{\langle e^2\cdot \rm{I}_2\rangle}$.

\section{Vaisman manifolds}\label{Vaisman_manifolds}
We  outline the basic facts needed for Vaisman manifolds. For a more detailed exposure we refer to \cite{Dr-Or} and recent papers by Verbitsky and the first author.

Let $(M,g,J)$ be a connected, smooth, Hermitian manifold of real dimension $2n+2$, with $n\ge 1$. Let $\omega(X,Y):=g(X,JY)$, with $X, Y \in \Gamma(TM)$, denote the fundamental 2-form.  We also denote by $\nabla$ the Levi-Civita connection of $\omega$.

\begin{definition}
	A Hermitian manifold $(M,J,g)$ is called {\bf Vaisman} if there exists a $\nabla$-parallel 1-form $\theta$ satisfying $d\omega=\theta\wedge\omega$.
\end{definition}

\begin{remark}\label{K_covering}
	The universal Riemannian cover of a Vaisman manifold is a metric cone of a Sasaki manifold (see e.g. \cite{Or-Ve}). Hence the local structure of a Vaisman manifold implies the existence of a local Sasaki structure transverse to the flow generated by the unitary vector field $U$.
\end{remark}

\example All
diagonal Hopf manifolds are Vaisman (\cite{ov_pams}). The Vaisman compact complex
surfaces are classified in \cite{Be}, see also \cite{_ovv:surf_}.

\remark There exist compact LCK manifolds which do not admit Vaisman metrics. Such are the LCK Inoue surfaces, \cite{Be}, the Oeljeklaus-Toma manifolds, \cite{oti2}, and the non-diagonal Hopf manifolds, \cite{ov_pams}, \cite{_ovv:surf_}.

\begin{remark}
As $\nabla\theta=0$ implies $d\theta=0$, a Vaisman manifold is locally conformally K\"ahler (LCK) \cite{Dr-Or}. The 1-forms $\theta$ is called the {\bf Lee form}.
\end{remark}

\medskip

The following is a fundamental observation due to I. Vaisman.

\begin{proposition} {\rm (\cite{Va})}\label{dthetac}
	The following  equation holds on a Vaisman manifold:
	\begin{equation} \label{theta_theta_c}
		\omega=d \theta^c- \theta\wedge \theta^c.
	\end{equation}
\end{proposition}

\subsection{The canonical foliation of a Vaisman manifolds}\label{canonical_foliation}

Since the Lee form is parallel, it has constant norm and hence we can normalize it such that  $\left\Vert \theta \right\Vert =1$. We denote by $\theta ^{c}=\theta \circ J$ the {\bf anti-Lee form}. Let $U=\theta ^{\sharp}$ and $V=\theta ^{c\sharp}$ denote the respective $g$-dual vector fields:
\[
\theta(U)=1,\quad \theta^c(V)=1,\quad \theta(V)=0,\quad \theta^c(U)=0.
\]

Note that the parallelism of $\theta$ implies \begin{equation}\label{conex_Vaisman}
	\nabla _UU=\nabla _VU=0,\qquad \nabla _UV=\nabla _VV=0.
\end{equation}

\remark\label{_canon_foli_totally_geodesic_Remark_}
On a Vaisman manifold, the Lee and anti-Lee fields  are real holomorphic ($\Lie_{U}I=\Lie_{V}I=0$) and Killing ($\Lie_{U}g=\Lie_{V}g=0$),\footnote{We denote with $\Lie_X$ the Lie derivative along the vector field $X$.} see \cite{Dr-Or}. Moreover, they commute: $[U,V]=0$. The 2-dimensional foliation $\Sigma$   they generate is called  {\bf the canonical foliation}. Its leaves are totally geodesic.

\remark
Recall that a metric defined on a foliated manifold is said to be {\bf bundle-like} if locally it can be identified with a Riemannian submersion \cite{Re}. The Vaisman metric is thus bundle-like  with respect to the canonical foliation \cite{Va}.

\medskip

Consider the canonical splitting of the tangent bundle
\begin{equation}\label{split_tang}
TM=Q\oplus \Sigma=Q\oplus \langle U\rangle\oplus \langle V\rangle,
\end{equation}
where $Q$ is the transverse $g$-orthogonal complement of $\Sigma$. Let  $\pi_Q:TM \rightarrow Q$ be  the canonical projection  induced by the splitting. This also imply a splitting of the metric
\[
g=g^T\oplus g^\Sigma,
\]
where $g^T$ is the transverse metric, and $g^\Sigma$ is the leafwise metric. By \cite{Va}, $g^\Sigma$ is Euclidean and $g^T$ the transverse metric, as in the case of a Riemannian submersion, can be locally projected on a local transversal (i.e.  a local submanifold of dimension $2n$, transverse to the leaves of the foliation). Moreover, $g^T$ is K\"ahler and hence $\Sigma$ is {\bf transversally K\"ahler}.

Along with the metric, we define other geometric objects that can be locally projected on submanifolds transverse to the leaves.

\definition
The de\thinspace Rham complex of the basic (projectable) differential
forms is defined as the restriction
\[
\Omega _b\left( M\right) :=\left\{ \alpha \in \Omega \left(
M\right) \mid \iota _X\alpha =0,\ \Lie_X\alpha =0 \ \ \mbox{\,for any\,}\ X\in \Gamma \left(\Sigma \right) \right\}.
\]
Here $\iota $ stands for interior product.

\definition
The {\bf basic de\thinspace Rham derivative} is the restriction of the
exterior derivative $d$, namely $d_b:=d_{\mid \Omega _b\left( M\right)}$
(see e.g.  \cite{Ton}).

\subsection{The transverse geometry of a Vaisman manifold}

On a Vaisman manifold the complex structure $J$ and the transverse fundamental 2-form $\omega^T$ are projectable. In particular, $\omega^T$ is a basic form, according to the above definition.

Since $\Sigma$ is bundle-like, we  choose an orthonormal frame  $\{e_{1},\ldots,e_n, e_{n+1,}\ldots,e_{2n},U,V\}$  such that
\[
\begin{split}
&J(U)=V,\quad\, \qquad J(V)=-U, \\
&J(e_j)=e_{n+j}, \quad J(e_{n+j})=-e_j.
\end{split}
\]
The dual frame will be denoted $\{e_{1,}^{\flat },\ldots,e_n^{\flat},e_{n+1,}^{\flat}, \ldots,e_{2n}^{\flat },\theta, \theta^c\dot \}$.

\subsubsection{Foliated charts on Vaisman manifolds}\label{foliated_map}

Consider a local foliated chart  $\left(\mathcal{U}, z^j, \bar{z}^j,x,y\right)$. Here  $\mathcal{U}=\mathcal{T}\times \mathcal{O}$, where:
\begin{itemize}
	\item   the local transversal $\mathcal{T}$ is a local chart  $\left(\mathcal{T}, z^j, \bar{z}^j\right)$ on a complex {\K} manifold,
$z^j$, $\bar{z}^j$ are the transverse complex coordinates, $z^j=x^j+\ii y^j$, $\bar{z}^j=x^j-\ii y^j$,
\item  $x,y$ are the leafwise real coordinates defined by $\frac \6 {\6x}=U,\  \frac \6 {\6y}=V$.
\end{itemize}
\begin{remark}\label{basic_vector}
The differential forms $dz^j$, $d\bar z^j$ are basic. Recall that $X\in \Gamma (Q)$ is {\bf a basic vector field} if $[X, Y]\in \Gamma (Q)$, for any leafwise tangent vector field $Y \in \Gamma(\Sigma)$ (see \eg \cite[Chapter 1]{Mo}). It follows that the vector fields $\frac \6 {\6z^i}$, $\frac \6 {\6\bar z^i}$ dual to $dz^j$, $d\bar z^j$ are basic.
\end{remark}

\subsubsection{The complex structure in a foliated chart}

The action of the transversal complex structure $J_0$ on the coordinate vector fields is:
In the local map, the standard complex structure $J_0$ is defined by the relations
\begin{equation*}
	\begin{split}
		J_0\left( \frac\partial{\partial x^j} \right)&=\frac\partial{\partial y^j}, \qquad
		J_0\left( \frac\partial{\partial y^j} \right)=-\frac\partial{\partial x^j}, \\
		J_0\left( \frac\partial{\partial x} \right)&=\frac\partial{\partial y}, \qquad
		J_0\left( \frac\partial{\partial y} \right)=-\frac\partial{\partial x}.
	\end{split}
\end{equation*}

As pointed out in the \ref{K_covering},  a Vaisman manifold is locally a Riemannian submersion, the fibers having a local structure of a Sasaki manifold (see e.g. \cite{Or-Ve}). The local structure of a Sasaki manifold using a local transverse {\K} potential is described in \cite{Go-Ko-Nu}  (see also \cite{Sm-Wa-Zh}). We consider a local K\"ahler potential $h$ (associated to the transverse K\"ahler structure). The local function $h$ is basic, constant along the leaves of the canonical foliation. Since the transverse {\K} structure of the Vaisman manifold is locally the same as the local {\K} structure of the Sasakian structure of the fibers, the complex structure $J$ is expressed in terms of the potential as:
\[
\begin{split}
  J&=J_0+\frac\partial{\partial x}\otimes d^c \left(h\right)+\frac\partial{\partial y}\otimes d^c \left(h\right)\circ J_0 \\
   &=J_0+\frac\partial{\partial x}\otimes \left(\ii \frac{\6 h}{\6 \bar{z}^j}d\bar{z}^j-\ii \frac{\6 h}{\6 z^j}dz^j\right)+\frac\partial{\partial y}\otimes \left(\ii \frac{\6 h}{\6 \bar{z}^j}d\bar{z}^j-\ii \frac{\6 h}{\6 z^j}dz^j\right)\circ J_0.
\end{split}
\]
In the coordinates of the foliated chart,
the Lee and anti-Lee forms are:
\[
\theta=dx, \quad \theta^c =-\theta \circ J=dx+ \ii \frac{\6 h}{\6z^j}dz^j-\ii \frac{\6 h}{\6\bar{z}^j}d\bar{z}^j.
\]
We introduce a new basis of local vector fields $X_1, \ldots  X_n, \bar{X}_1, \ldots  \bar{X}_n$ in $Q$, in which the complex structure acquire a simpler expression. Let
\begin{equation}\label{X_j}
 X_j =\frac{\6}{\6 z^j}-\ii \frac{\6 h}{\6 z^j}\frac{\6}{\6 x}, \qquad
 \bar{X}_j =\frac{\6}{\6 z^j}+\ii \frac{\6 h}{\6 \bar{z}^j}\frac{\6}{\6 x}.
\end{equation}
A direct computation shows that
\begin{equation}\label{J_X}
\begin{split}
J(X_j)=\ii X_j,\quad J(\bar X_j)=-\ii \bar X_j.
\end{split}
\end{equation}
We then work with the local basis of vector fields $\{X_j,\bar{X}_j, U, V\}$. The dual basis of 1-forms is $\{dz^j, d\bar{z}^j, \theta, \theta^c\}$.

\begin{remark}
The local vector fields $X_j,\bar{X}_j$ and the dual one-forms $dz^j, d\bar{z}^j$ are projectable on a local transversal. 
\end{remark}

\medskip

Let $g^T_{j\bar{k}}:=g(X_j, \bar{X}_k)$ be the coefficients of the transverse metric in these new coordinates. We then have:
\begin{equation*}
    g^T=g^T_{j\bar k}dz^j \odot d\bar z^k:=g^T_{j\bar k}(dz^j \otimes d\bar z^k
    +d\bar z^k \otimes dz^j).
\end{equation*}
The associated  transverse fundamental 2-form is $\omega^T=\ii g^T_{j\bar k}dz^j \wedge \bar z^k$.

\section{Deformation of a Vaisman structure using a differential basic 1-form}\label{deform}

\subsection{Description of the method}\label{method}
Let $(M,J,g, \theta)$ be a closed Vaisman manifold, with real dimension $2n+2$. In the sequel we construct a deformation $(M,J_{t},g_{t})$ of the Vaisman structure for $t\in (-\varepsilon ,\varepsilon )$, with $J_{0}=J$, $g_{0}=g$, in such a way that the canonical foliation remain unchanged and the deformation concerns only the transverse orthogonal complement and the transverse K\"ahler geometry.

Let $\{\zeta _{t}\}_{t}$ be a family of differential 1-forms, which are basic with respect to $\Sigma$  ($\zeta _{t}(U)=\zeta _{t}(V)=0$, $\Lie_{U}\zeta _{t}=\Lie_{V}\zeta _{t}=0$), with  $\zeta _{0}=0$.
Later on, we shall add new conditions on the 1-forms $\{\zeta _{t}\}_{t}$.

We fix $U_{t}=U$, $V_{t}=V$. This ensures that $\Sigma$ is not affected by the deformation.

Consider the canonical splitting \eqref{split_tang} of the tangent bundle. Now  deform the normal bundle $Q$ to
\begin{equation*}
Q_{t}=\left\{ v-\zeta _{t}(v)V\mid v\in Q\right\} .
\end{equation*}
and associate  the morphism of vector (sub-)bundles
$$\pi _{t}:Q\rightarrow Q_{t},\quad \pi _{t}(v)=v-\zeta _{t}(v)V, \quad  \,v\in Q.$$
Then
$$\pi_{t}^{-1}(w)=w+\zeta _{t}(w)V, \quad w\in Q_{t}.$$
Define $J_t$ on $\Sigma$ by
\[
J_{t}(U)=V,\,J_{t}(V)=-U.
\]
and on $Q_t$ by $J_{t}|_{Q_t}=\pi _{t}\circ J\circ \pi _{t}^{-1}$. Then
\begin{equation}
\begin{split} \nonumber
J_{t}|_{Q_t}(w) &=\pi _{t}\circ J\circ \pi _{t}^{-1}(w) \\
&=\pi _{t}\circ J(w+\zeta _{t}(w)V) \\
&=J(w+\zeta _{t}(w)V)-\zeta _{t}(J(w+\zeta _{t}(w)V)) \\
&=J(w)-\zeta _{t}(w)U-\zeta _{t}(J(w))V,
\end{split}
\end{equation}
Written in a compact way, the deformation $J_t$ reads
\begin{equation}
	J_{t}=J-U\otimes \zeta _{t}-V\otimes (\zeta _{t}\circ J).  \label{J_t}
\end{equation}
\begin{claim}
$J_{t}^{2}=J^{2}=-\mathrm{Id}$. Hence $J_{t}$ is an almost complex structure for all $t\in\R$.
\end{claim}

\proof
Since $\zeta _{t}$ is a basic  1-form, we have
\begin{equation*}
\begin{split}
J_t^2&=\left(J-U\otimes \zeta _{t}-V\otimes (\zeta _{t}\circ J) \right)^2 \\
&=J^2 -J(U)\otimes\zeta_t- J(V)\otimes\zeta_t\circ J-U\otimes\zeta_t \circ J+V\otimes\zeta_t \\
&=J^2=-\Id.
\end{split}
\end{equation*}
\hfill{\endproof}

\begin{proposition}\label{J_integrabil}
If $d\zeta_t$ are of type (1,1) with respect to $J$, then $J_t$ is integrable for all $t\in\mathbb{R}$.
\end{proposition}
\begin{proof} Let $N^{J_t}$ be the Nijenhuis tensor field associated to $J_t$:
	\begin{equation*}
		N^{J_t}(\cdot ,\cdot ) =J_{t}^{2}[\cdot ,\cdot ]-J_{t}[J_{t}\cdot ,\cdot]-J_{t}[\cdot ,J_{t}\cdot ] +[J_{t}\cdot ,J_{t}\cdot ]
	\end{equation*}
We show that $N^{J_t}$ vanishes. It is enough to show that it vanishes on basic vector fields. Let then  $X$ and  $Y$ be basic vector fields (we include here the case when $X$, $Y$ are tangent to the leaves). Then $[U,X]$, $[U,Y]$, $[V,X]$ and $[U,Y]$ are also tangent to the leaves of the foliation. As $\zeta _{t}$ is basic, we have
\begin{equation}
\begin{split}
N^{J_t}(X,Y) &=N^J(X,Y)-d\zeta _{t}(JX,Y)U-d\zeta _{t}(X,JY)U \\
&+d\zeta _{t}(X,Y)V-d\zeta _{t}(J(X),J(Y))V \\
&-\zeta _{t}(X)(\Lie_{U}J)(Y)-\zeta _{t}(J(X))(\Lie_{V}J)(Y) \\
&+\zeta _{t}(J(Y))(\Lie_{V}J)(X)+\zeta _{t}(J(Y))(\Lie_{U}J)(X) \\
&-\zeta _{t}(N^J(X,Y)).
\end{split}
\end{equation}

Recall the relations
\begin{equation}\label{J_Lie_inv}
\Lie_{U}J=\Lie_{V}J=0.
\end{equation}
Since $d\zeta_t$ is of type (1,1) and $J$ is integrable, by  \eqref{J_Lie_inv} we obtain
\begin{equation*}
N^{J_t}(X,Y)=0,
\end{equation*}
for $X$ and $Y$ basic vector fields.
\hfill\end{proof}

\hfill

\fbox{From now on, we assume that $d\zeta_t$ are of type (1,1) w.r.t $J$.}

\hfill

We now choose $\theta _{t}=\theta $, so $\theta _{t}^{c}=\theta \circ
J_{t}=\theta ^{c}+\zeta _{t}$ is the deformation of anti-Lee 1-form. One can easily check
\begin{equation*}
\theta _{t}^{c}(V)=1,\quad \theta _{t}^{c}(U)=0.
\end{equation*}
\begin{remark}
Recall (\ref{dthetac}) that on a Vaisman manifold,
\begin{equation}\label{dtetac}
	d\theta^c=\omega+\theta\wedge\theta^c.
\end{equation}
Therefore $(d\theta^c)^{n}\wedge\theta\wedge\theta^c$ is a volume form (where $\dim_\R M=2n+2$). We need that  the form $(d\theta _{t}^{c})^{n}\wedge 	 \theta_{t}^{c}\wedge \theta $ remain a volume form during the deformation process (see also \cite[pag. 447]{Bo-Ga}). This is not automatic, hence we  add a condition on $\{\zeta _{t}\}_{t}$ such as to fulfill this requirement.

We have $\|\theta\|=\|\theta^c\|=\|(d\theta^c)^{n}\wedge\theta\wedge\theta^c\|=1$, $\|d\theta^c\|=\sqrt{n}$ w.r.t. the initial metric $g$. In the expansion

\[
(d\theta _{t}^{c})^{n}\wedge \theta\wedge \theta_{t}^{c}=
\sum^n_{k=1} \binom nk (d\theta _{t}^{c})^{n-k}\wedge d\zeta^k_t \wedge\theta \wedge \theta^c
+\sum^n_{k=1} \binom nk (d\theta _{t}^{c})^{n-k}\wedge d\zeta^k_t \wedge \theta \wedge \zeta_t
\]
there are $2^{n+1}$ terms in the right hand side. If we choose $\zeta_t$ such that
\begin{equation}\label{estimare1}
\|\zeta_t\|<\frac 1{(2^{n+1}-1){n}^\frac{n-1}{2}}, \quad\|d\zeta_t\|<\frac 1{(2^{n+1}-1){n}^\frac{n-1}{2}},
\end{equation}
then, using the Cauchy-Buniakovski-Schwarz inequality $\|\alpha^1\wedge \alpha^1\|\le \|\alpha^1\|\cdot|\alpha^2\|$ for differential forms $\alpha^1$, $\alpha^2$, we get the estimate
\[
\begin{split}
\|(d\theta _{t}^{c})^{n}\wedge \theta\wedge \theta_{t}^{c}\|&\ge 1-
\sum^n_{k=2} \binom nk \|d\theta _{t}^{c}\|^{n-k}\cdot\|d\zeta^k_t\| \cdot \|\theta \| \cdot \|\theta^c\| \\
&-\sum^n_{k=1} \binom nk \|d\theta _{t}^{c}\|^{n-k}\cdot \|d\zeta^k_t\| \cdot \|\theta \| \cdot \|\zeta_t\| \\
&> 1-(2^{n+1}-1)\frac1 {2^{n+1}-1}> 0.
\end{split}
\]
Then $(d\theta _{t}^{c})^{n}\wedge \theta\wedge \theta_{t}^{c}$ remains a volume form.

In conclusion, to assure \eqref{estimare1} it is enough to divide  any given $\zeta_t$  by a constant greater than $(2^{n+1}-1){n}^\frac{n-1}{2}\max\{\|\zeta_t\|, \|d\zeta_t\|\}$.

\end{remark}

\hfill

\fbox{In the sequel we assume that $\zeta_t$ fulfill the inequalities \eqref{estimare1}}.

\hfill

\begin{claim} \label{J_t_J}
	 Let $\omega_t$ be the differential 2-form
\begin{equation*}
\omega _{t} =d\theta _{t}^{c}-\theta \wedge \theta _{t}^{c} =d\theta ^{c}+d\zeta _{t}-\theta \wedge \theta _{t}^{c}.
\end{equation*}
Then  $\omega_t$ is of type (1,1) w.r.t. $J_t$.
\end{claim}
\begin{proof}
	It will be  enough to verify the claim on basic 1-forms.
	
{\noindent \bf Step 1: $J_t=J$ on basic differential forms.}  Let $\alpha$ be a basic 1-form. Then
\[
J_t(\al)(v)=-\al(Jv-\zeta_t(v)U-\zeta(Jv)V)=-\al(Jv)=J(\al)(v),
\]
for any $v\in \Gamma(TM)$, since $\al(U)=\al(V)=0$.

{\noindent \bf Step 2: $d\zeta_t$ and $d \theta^c$ are of type (1,1) w.r.t. $J_t$.} Indeed, since $\zeta_t$ are basic and  $[\Lie, d]=0$, we see that $d\zeta_t$ are  basic too. By assumption, $d\zeta_t$ are of type $(1,1)$ w.r.t. the complex structure $J$. Then by Step 1  they are of type (1,1) w.r.t. $J_t$, for all $t\in\RR$.

By \eqref{dtetac}, $d\theta^c$   is  basic  and of type $(1,1)$ w.r.t. $J$, and hence $\omega_t$ is a sum of (1,1) 2-forms.
\hfill\end{proof}

\smallskip

Now, since $\theta $ is closed, we have
\begin{equation}
	d\omega _{t}=\theta \wedge d\theta _{t}^{c}=\theta \wedge \omega _{t}.
	\label{omega_t}
\end{equation}
and hence
\begin{claim}
The manifold $(M,\omega _{t}, \theta)$ is locally conformally symplectic  for all $t\in\R$. \hfill{\endproof}
\end{claim}

\medskip

We define a symmetric (0,2) tensor field by
\begin{equation}\label{g_t1}
\begin{split}
g_{t} &=\omega _{t}\circ (J_{t}\otimes \mathrm{Id}) \\
&=d\theta ^{c}\circ (J_{t}\otimes \mathrm{Id})+d\zeta _{t}\circ (J_{t}\otimes \mathrm{Id}) +\theta _{t}^{c}\otimes \theta _{t}^{c}+\theta \otimes \theta \\
&=d\theta ^{c}\circ (J\otimes \mathrm{Id})+d\zeta _{t}\circ (J\otimes \mathrm{Id}) +\theta _{t}^{c}\otimes \theta _{t}^{c}+\theta \otimes \theta.
\end{split}
\end{equation}
In the last equality we used the fact that $d\theta ^{c}$, $d\zeta _{t}$ are basic differential forms.
Note that
\begin{equation} \nonumber
\theta _{t}^{c}\otimes \theta _{t}^{c} =\theta ^{c}\otimes \theta^{c}
+\zeta _{t}\otimes \theta ^{c}+\theta ^{c}\otimes \zeta _{t}
+\zeta _{t}\otimes \zeta _{t}.
\end{equation}
The tensor field $g_t$ is related  to the metric tensor $g$ by the formula:
\begin{equation}
g_{t}=g+d\zeta _{t}\circ (J\otimes \mathrm{Id})  \label{g_t}
+\zeta _{t}\otimes \theta ^{c}+\theta ^{c}\otimes \zeta _{t}
+\zeta _{t}\otimes \zeta _{t}.
\end{equation}

\begin{remark}
We want $g_t$ to be a metric such as to regard \eqref{g_t1} as a deformation of the initial metric $g$. To this end, we need to find a condition which implies the positivity of the tensor field $g_t$. At any point $x \in M$ consider $S_x^1M:=\{v\in T_xM \mid \|v\|=1\}$, where the norm is taken w.r.t. the metric $g$. Define $$\mu_x:=\min_{v\in S^1_xM} d\zeta(Jv,v).$$ As $S^1_xM$ is a compact set, $\mu_x>-\infty$. Moreover, $\mu_x$ is a smooth function on the compact manifold $M$, and attains its minimum at a certain point. We make the assumption that

\begin{equation}\label{estimare2}
\min_{x\in M} \mu_x>-1.
\end{equation}

Note that by \eqref{g_t1}, $g_t(v,v)>0$ for all $v$  tangent to the leaves of the canonical foliation.

For a transverse unit vector field $v\in Q_{t,x}$, we have
\[
g_t(v,v)=d\theta^c(Jv,v)+d\zeta(Jv,v)=1+d\zeta(Jv,v)>0
\]
by \eqref{estimare2} and since $d\theta^c(Jv,v)$ is the transverse part of the metric $g$. Hence \eqref{estimare2} is a sufficient condition for $g_t$ to be positive definite.
\end{remark}

\begin{remark}
	Let $c>|\min_{x\in M} \mu_x|$. Then, for any  given $\zeta_t$, the form $c\cdot\zeta$ satisfies the inequality \eqref{estimare2}.
\end{remark}

\hfill

\fbox{In the following we assume that $\zeta_t$ satisfy the inequality \eqref{estimare2}.}

\hfill

The form $\zeta _{t}$ is basic, and hence $\Lie_{U}\zeta _{t}=0$, $\Lie_{U}d\zeta _{t}=0$; then from \eqref{g_t} we obtain
\begin{equation}
\Lie_{U}g_{t}=0,  \label{U_Killing}
\end{equation}
so $U$ stays Killing w.r.t. all the metrics $g_t$.

A straightforward computation shows  that $U$ and $V$ are vector fields metrically equivalent to $\theta $ and $\theta _{t}^{c}$, with respect to the metric $g_{t}$. Also,
\begin{equation}
\left\Vert \theta \right\Vert _{g_t}=\left\Vert \theta _{t}^{c}\right\Vert_{g_t}=1.  \label{norma_const}
\end{equation}
Gathering together the above facts we obtain:

\begin{theorem}\label{g_t_vai}
	Let $(M,J,g,\theta)$ be a Vaisman manifold and $J_t, g_t$ given by \eqref{J_t} and \eqref{g_t}. Assume that the basic 1-forms $\zeta_t$ satisfy \eqref{estimare1} and \eqref{estimare2} and  $d\zeta_t$ are of type (1,1) w.r.t. $J$. Then  $(M,
J_{t}, g_{t})$ is a Vaisman structure with Lee form $\theta$ for all $t\in\R$.
\end{theorem}

\subsection{Construction of the basic forms $\zeta_t$} \label{Construction}

We start with a globally defined, basic function $\f $ (w.r.t. the canonical foliation). Let  $d_b^c:=Jd_bJ^{-1}$ be the corresponding complex operator acting on basic $k$-forms. By Step 1 in the proof of \ref{J_t_J},  $d_{b}^{c}$ remains unchanged by the deformation of the complex structure, so this operator can be taken w.r.t. the complex structure $J$.Let
$$\zeta _{t}:=td_{b}^{c}\f.$$
Clearly $td_{b}^{c}\f$ is a basic 1-form of type $(1,1)$, and we take $t$ small enough such that $td_{b}^{c}\f$ satisfies the inequalities \eqref{estimare1} and \eqref{estimare2}.
Then, according to \ref{g_t_vai}, the family of metrics $g_t$ defined by
\begin{equation}
g_{t}=g+t(dd_{b}^{c}\f\circ (J\otimes \mathrm{Id})  \label{g_tt}
+d_{b}^{c}\f\otimes \theta ^{c}+\theta ^{c}\otimes d_{b}^{c}\f)
+t^2d_{b}^{c}\f\otimes d_{b}^{c}\f.
\end{equation}
represent a deformation of the initial Vaisman metric $g$.

\subsection{Example: deforming a diagonal Hopf surface}\label{ex1}
We describe the geometry of a classical Hopf surface as it is useful for our deformation process. For the reader's convenience, we provide the explicit computations.

Let $\Delta $ the cyclic group  generated by the holomorphic transformation $z\mapsto e^2z$ on $\CC^2\setminus\{0\}$. The quotient space $H:=\left( \CC^2\setminus\{0\}\right) /\Delta $, is a (primary) {\bf Hopf surface}.

To describe the  foliated coordinates, consider the Hopf fibration $S^3 \rightarrow \CC P^1$ and let $w:=u+iv$ be the complex local coordinate on $S^2\simeq \CC P^1$. The standard {\K} metric on  the space of leaves $S^2$ is determined by the local potential $\frac 12\log (1+\left|w \right|^2) $.  Let $x$ be the coordinate on the factor $ S^{1}$ and $y$ the coordinate of the fibre of the Hopf fibration. Then $\left(w, x,y \right)$ are local coordinates on the foliated manifold (with $w$ the transverse coordinate). 

In the local foliated chart, the standard complex structure $J_0$ is defined by the relations
\begin{equation*}
\begin{split}
J_0\left( \frac\partial{\partial u} \right)&=\frac\partial{\partial v}, \qquad
J_0\left( \frac\partial{\partial v} \right)=-\frac\partial{\partial u}, \\
J_0\left( \frac\partial{\partial x} \right)&=\frac\partial{\partial y}, \qquad
J_0\left( \frac\partial{\partial y} \right)=-\frac\partial{\partial x}.
\end{split}
\end{equation*}

Define the almost complex structure $J$ on $H$ as
\[
  J=J_0+\frac\partial{\partial x}\otimes d^c \left(\frac1 2 \log (1+\left|w \right|^2) \right)+\frac\partial{\partial y}\otimes d^c \left(\frac1 2 \log (1+\left|w \right|^2) \right)\circ J_0.
\]
By \ref{J_integrabil}, $J$ is integrable. Since $\left[\frac {\6}{\6 x}, \frac {\6}{\6 y}\right]=0$ and the {\K} potential is basic, $J$ is invariant along the leaves:
\begin{equation}\label{Lie_invar}
\Lie_{\frac {\6}{\6 x}}J=\Lie_{\frac {\6}{\6 y}}J=0.
\end{equation}
In these  local coordinates, the Lee form is $\theta=dx,$ and anti-Lee is
\[
\begin{split}
\theta^c&=-\theta\circ J \\
&= -dx \left( J_0+\frac\partial{\partial x}\otimes d^c \left(\frac1 2 \log (1+\left|w \right|^2) \right) +\frac\partial{\partial y}\otimes d^c \left(\frac1 2 \log (1+\left|w \right|^2) \right)\circ J_0 \right) \\
&= -dx\circ J_0 +d^c \left(\frac1 2 \log (1+\left|w \right|^2)  \right) \\
&= dy-\frac \ii 2 \frac w {1+\left|w \right|^2} d\bar{w} + \frac \ii 2 \frac {\bar{w}} {1+\left|w \right|^2} dw.
\end{split}
\]
Then the fundamental 2-form reads
\[
\begin{split}
\omega=&d\theta^c -\theta\wedge \theta^c \\
=&-d d^c \left(\frac1 2 \log (1+\left|w \right|^2 \right) - dx \wedge \left(dy-\frac \ii 2 \frac w {1+\left|w \right|^2} d\bar{w} + \frac \ii 2 \frac {\bar{w}} {1+\left|w \right|^2} dw \right) \\
=& -\ii\frac 1 {(1+\left|w \right|^2)^2} dw\wedge d\bar w + dx\wedge dy +\frac{\ii}{2} \frac w {1+\left|w \right|^2} dx\wedge d\bar{w} \\
&-\frac{\ii}{2} \frac {\bar{w}} {1+\left|w \right|^2} dx\wedge dw.
\end{split}
\]
The corresponding metric is $g=\omega\circ (J\otimes \mathrm{Id})$. It is positive definite since the transverse part is the standard {\K} metric on $S^2$ while on leaves it is the Euclidean metric.  Using \eqref{Lie_invar} one easily proves that 
\[
\Lie_{\frac {\6}{\6 x}}g=\Lie_{\frac {\6}{\6 y}}g=0,
\]
and hence the metric is bundle-like with respect to the foliated structure. Since
\[
\left\|\frac{\6}{\6 x}\right\|=\left\|\frac{\6}{\6 y}\right\|=1,
\]
the metric $g$ is Vaisman.

We deform the Vaisman structure, using a smooth basic function defined on our foliation, as indicated in Subsection \ref{Construction}. Let
\[
f(w,x,y)=\frac 1 2 \left(\frac {1-\left|w \right|^2}{1+\left|w \right|^2} \right)^2,
\]
and take $\zeta_t:=t\cdot f$.
\begin{remark}
On the manifold $S^3$,  consider the standard  parametrization
\[
 \left[ z^1,z^2\right]=\left[ \cos s\, e^{i\Phi_1}, \sin s\, e^{i\Phi_2} \right].
\]
One can see that $f=\frac 1 2 \cos ^2\! s$ is a basic, globally defined function on  $S^3$ (see \textrm{e.g.} \cite[Example 1]{Slo}); the same holds on the Hopf manifold $H$.
\end{remark}

\medskip

The family of complex structures constructed as in Subsection \ref{method} is the following:
\[
\begin{split}
  J_t=&J-\frac\partial{\partial x}\otimes d^c \left(\frac1 2 t^2 \left(\frac {1-\left|w \right|^2}{1+\left|w \right|^2}\right)^2 \right) -\frac\partial{\partial y}\otimes d^c \left(\frac1 2 t^2 \left(\frac {1-\left|w \right|^2}{1+\left|w \right|^2}\right)^2 \right)\circ J.
\end{split}
\]
The Lee form  remains fixed, $\theta_t=\theta=dx$, and the anti-Lee 1-form  is
\[
\begin{split}
\theta_t^c&=-\theta\circ J_t \\
&=- dx \left( J-\frac\partial{\partial x}\otimes d^c \left(\frac1 2 t \left(\frac {1-\left|w \right|^2}{1+\left|w \right|^2}\right)^2 \right) -\frac\partial{\partial y}\otimes d^c \left(\frac1 2 t \left(\frac {1-\left|w \right|^2}{1+\left|w \right|^2}\right)^2 \right)\circ J \right) \\
&= -dx\circ J +d^c \left(\frac1 2 \log (1+\left|w \right|^2)  \right) \\
&= dy-\frac \ii 2 \left(\frac w {1+\left|w \right|^2}-t \frac{1-\left|w\right|^2}{(1+\left|w\right|^2)^3}w\right) d\bar{w}
+\frac \ii 2 \left(\frac {\bar{w}} {1+\left|w \right|^2}-t \frac{1-\left|w\right|^2}{(1+\left|w\right|^2)^3}\bar{w}\right) dw.
\end{split}
\]
Finally, we describe the deformation of the fundamental 2-form:
\[
\begin{split}
\omega_t=&d\theta^c_t -\theta\wedge \theta_t^c \\
=& -\ii\frac{\left[(1-t)+(2+4t)\left\vert w\right\vert ^{2}+(1-t)\left\vert
w\right\vert ^{4}\right]^2 }{\left( 1+\left\vert w\right\vert ^{2}\right) ^{8}} dw\wedge d\bar w - dx\wedge dy  \\
&+\frac{\ii w(1+\|w\|^2-t(1-\|w\|^2))}{2(1+\|w\|^2)^3} dx\wedge d\bar{w} \\
&-\frac{\ii\bar{w}(1+\|w\|^2-t(1-\|w\|^2))}{2(1+\|w\|^2)^3}dx\wedge dw.
\end{split}
\]

If $M_{t}$ is the skew-symmetric matrix associated to $\omega _{t}$ with
respect to the local basis $(du,dv,dx,dy)$, then
\[
\det M_{t} =\frac{\left[(1-t)+(2+4t)\left\vert w\right\vert ^{2}+(1-t)\left\vert
w\right\vert ^{4}\right]^2 }{\left( 1+\left\vert w\right\vert ^{2}\right) ^{8}}.
\]

Taking $t\in(-\frac 1 2,\frac 1 2)$ we see that all terms that appear are strictly positive and $\det M_{t}\neq 0$  for all $w$. Hence  $\omega _{t}$ is non-degenerate for all $t\in (-1/2, 1/2)$. Define the family of metrics  $g_{t} =\omega _{t}\circ (J_{t}\otimes \mathrm{Id})$.

\begin{claim}$g_t$ is positive definite.
\end{claim}
\begin{proof}
Denote by $\{\lambda_i(t,w,x,y)\}_{1 \le i \le n}$ the eigenvalues associated to $g_t$ at the point of coordinates $(w,x,y)$. As $g_t$ is not degenerate, the smooth eigenfunctions $\lambda_i$ are nowhere vanishing, so they have constant sign. Moreover, $\lambda_i>0$ since $\lambda_i(0,w,x,y)>0$, hence the claim is proved.
\hfill\end{proof}

\hfill

{\small
	
	\noindent {\sc Liviu Ornea\\
		University of Bucharest, Faculty of Mathematics and Informatics, \\14
		Academiei str., 70109 Bucharest, Romania}, and:\\
	{\sc Institute of Mathematics ``Simion Stoilow" of the Romanian
		Academy,\\
		21, Calea Grivitei Str.
		010702-Bucharest, Romania\\
		\tt lornea@fmi.unibuc.ro,   liviu.ornea@imar.ro}
	
	\hfill

	\noindent {\sc Vladimir Slesar\\
		 University ``Politehnica'' of Bucharest, Faculty of Applied Sciences \\
			313 Splaiul Independenţei,
			060042 Bucharest, Romania}, and\\
		{\sc University of Bucharest, Faculty of Mathematics and Informatics, \\14
			Academiei str., 70109 Bucharest, Romania}\\
	\tt vladimir.slesar@upb.ro }

\end{document}